# Termination of Picard Iteration for Coupled Neutronics/Thermal-Hydraulics Simulations

Dean Wang[1] and Paul K. Romano[2]

[1]Nuclear Engineering Program, The Ohio State University, Columbus, OH 43210
[2]Argonne National Laboratory, Lemont, IL 60439
wang.12239@osu.edu, promano@anl.gov

## INTRODUCTION

In our previous work [1], we performed a formal Fourier analysis (FA) of Picard iteration for the coupled nonlinear neutronics/thermal hydraulics (N/TH) problem, and derived theoretical predictions for the spectral radius of Picard iteration, which is a function of various parameters such as the temperature difference between the fuel and the coolant, the temperature coefficients of cross sections (i.e., Doppler feedback), the scattering ratio, and the core height. We also derived an estimate for underrelaxation based on Fourier analysis. In the analysis, we need to assume that during each Picard iteration, the neutronics and thermal hydraulics solutions are computed exactly or fully converged if an iterative method is used. One potential drawback of Picard iteration is that we may unnecessarily over solve the sub physics problems while the global convergence has not been reached. Accordingly, this raises the question: How accurate do we need to solve the subproblems to obtain an optimal computational efficiency?

In the case of Newton's method for solving a nonlinear system, inexact Newton methods have been successfully developed to address the over-solving issue [2,3]. At each outer Newton step, the inner iteration for solving the linear equation system is terminated when a relative tolerance criterion in the linear residual is satisfied. Birken analyzed inexact fixed-point iterations (including Picard iteration) by modeling the iteration as a perturbed fixed-point iteration [4], and proved that if the iteration converges, it converges to the exact solution irrespective of the tolerance in the inner systems, provided that a nonstandard relative termination criterion is employed. Recently, Senecal and Ji proposed a modified Picard iteration coupling method with adaptive, inexact termination criteria for the underlying single-physics codes [5]. Within Picard Iteration, inexact Newton methods are applied in the single-physics codes.

In this paper, we consider the coupled N/TH problem, in which the termination criterion for the neutronics iteration adopts an adaptive tolerance with respect to the fuel temperature residual at each Picard iteration. We refer to this coupling scheme as the inexact Picard iteration method. Fourier analysis is performed to investigate how the convergence behavior of Picard iteration is influenced by the inexact neutronics solution. It is found that if the convergence of the inner neutronics iteration is slow, Picard coupling may become unstable unless a tighter tolerance is used for the neutronics iteration. Nevertheless, our analysis indicates that a certain amount of over-solving is necessary for maintaining the stability of Picard iteration if the iterative solution of the subproblem is not fast enough. However, this issue has not been addressed in the previous studies.

## FORMULATION AND ALGORITHM

As in our previous work [1], we consider the following simple one-group, planar-geometry k-eigenvalue problem on the domain $0 \leq x \leq L$ with reflective boundary conditions:

$$\mu \frac{\partial \psi(x,\mu)}{\partial x} + \Sigma_t(T)\psi(x,\mu)$$
$$= \frac{1}{2}\Sigma_s(T)\phi(x) + \frac{1}{2k_{eff}}\nu\Sigma_f(T)\phi(x) , \qquad (1)$$

and the simplified heat transfer equation for a single typical PWR fuel pin:

$$T = T_m + A\Sigma_f(T)\phi(x) , \qquad (2)$$

with

$$A = \pi r_{fo}^2 \kappa R_t , \qquad (3a)$$

and

$$R_t = \left[\frac{1}{8\pi k_f} + \frac{1}{2\pi r_g h_g} + \frac{1}{2\pi k_c}\ln\left(\frac{r_{co}}{r_{ci}}\right) + \frac{1}{2\pi r_{co} h}\right] . \qquad (3b)$$

Picard iteration for the coupled N/TH system is written as

$$\mu \frac{\partial \psi^{(k+1,n+1)}(x,\mu)}{\partial x} + \Sigma_t(T^{(k)})\psi^{(k+1,n+1)}$$
$$= \frac{1}{2}\Sigma_s(T^{(k)})\phi^{(k+1,n)}(x) + \frac{1}{2k_{eff}^{(k+1,n)}}\nu\Sigma_f(T^{(k)})\phi^{(k+1,n)} , \qquad (4)$$

$$T^{(k+1)} = T_m + A\Sigma_f(T^{(k)})\phi^{(k+1,n+1)}(x) , \qquad (5)$$

where the superscript $k$ denotes the global Picard iteration number, and $n$ is the index for the inner neutronics iteration.

The algorithm for the inexact Picard iteration is described as follows.

| Algorithm I: Inexact Picard Iteration |
| --- |
| while $r_T^{(k)} = \frac{\|(T^{(k)}-T^{(k-1)})\|}{\|(T^{(k)})\|} > \epsilon$ <br>    while $r_N^{(k+1,n)} = \frac{\|(\phi^{(k+1,n)}-\phi^{(k+1,n-1)})\|}{\|(\phi^{(k+1,n+1)})\|} > \tau \cdot r_T^{(k)}$ <br>      $\boldsymbol{M}(T^{(k)})\phi^{(k+1,n+1)} = \frac{1}{k_{eff}^{(k+1,n)}}\boldsymbol{F}(T^{(k)})\phi^{(k+1,n)}$ <br>    end <br>    $T^{(k+1)} = \boldsymbol{\mathcal{T}}\phi^{(k+1,n+1)}$ <br> end |

In this simple model problem, the neutronics solution is obtained using an iterative method such as the power iteration, while the TH solution is obtained exactly. Since the final solution is generally not available, in practical implementation Picard iteration is terminated by checking if the flux residual between the two successive iterations is smaller than the prescribed tolerance. $r_T^{(k)}$ is the relative fuel temperature residual at the $k$-th Picard iteration. $r_N^{(k+1,n)}$ is the relative flux residual at the inner $n$-th neutronics iteration. $M$ represents the neutron migration and loss operator, and $F$ is the fission operator. $\mathcal{T}$ denotes the TH solver operator. The stopping criterion for the inner iteration is set to be proportional to the outer fuel temperature residual by using a "forcing" parameter $\tau$. Such an adaptive approach leads to an automatic tightening of the tolerance over Picard iteration. The forcing parameter $\tau$ will be determined by Fourier analysis below.

## PERTURBED PICARD ITERATION

In Fourier analysis, the inexact Picard iteration is modeled as a perturbed fixed-point iteration [4,6]:

$$\mu \frac{\partial \psi^{(k+1)}(x,\mu)}{\partial x} + \Sigma_t(T^{(k)}) \psi^{(k+1)}$$
$$= \frac{1}{2}\Sigma_s(T^{(k)})\phi^{(k+1)}(x) + \frac{1}{2k_{eff}}\nu\Sigma_f(T^{(k)})\phi^{(k+1)}, \quad (6)$$

$$T^{(k+1)} = T_m + A\Sigma_f(T^{(k)})\phi^{(k+1)}(x)\left[1 + \varepsilon^{(k+1)}\right]. \quad (7)$$

In the inner neutronics calculation, the solution is now assumed to be fully converged. However, when the updated neutron flux is passed to the TH model, it contains a small perturbation to account for the error of partial convergence. The perturbation parameter, $\varepsilon^{(k+1)}$, shall be determined based on the adaptative termination criterion for the inner neutronics solution, which is described as follows.

In Algorithm I, the termination criterion for the neutronics iteration is

$$\frac{\|\phi^{(k+1,n+1)} - \phi^{(k+1,n)}\|}{\|\phi^{(k+1,n+1)}\|} \leq \tau \frac{\|T^{(k)} - T^{(k-1)}\|}{\|T^{(k)}\|}. \quad (8)$$

Assuming that the neutronics iteration has linear convergence (which is often the case), we can estimate the iterative error of the neutron flux as

$$\|\phi^{(k+1,n+1)} - \phi^{(k+1,n)}\| \approx \frac{1-\rho_N}{\rho_N}\|\phi^{(k+1,n+1)} - \phi^{(k+1)}\|, \quad (9)$$

where $\phi^{(k+1)}$ is the fully converged flux and $\rho_N$ is the spectral radius of the neutronics iteration. In fact, the above error estimate is very sharp.

Similarly, we can have the error estimate for the fuel temperature as

$$\|T^{(k)} - T^{(k-1)}\| \approx \frac{1-\rho}{\rho}\|T^{(k)} - T_0\|, \quad (10)$$

where $T_0$ is the converged fuel temperature and $\rho$ is the spectral radius of Picard iteration.

Combining Eqs. (8) and (9), we obtain

$$\|\phi^{(k+1,n+1)} - \phi^{(k+1)}\|$$
$$\leq \tau \frac{\rho_N}{1-\rho_N} \frac{1-\rho}{\rho} \frac{\|\phi^{(k+1,n+1)}\|}{\|T^{(k)}\|} \|T^{(k)} - T_0\|. \quad (11)$$

Thus, the estimate for the perturbation is

$$\varepsilon^{(k+1)} = \pm\tau \frac{\rho_N}{1-\rho_N} \frac{1-\rho}{\rho} \frac{\|\phi^{(k+1,n+1)}\|}{\|\phi^{(k+1)}\|} \frac{\|T^{(k)} - T_0\|}{\|T^{(k)}\|}$$
$$\approx \pm\tau \frac{\rho_N}{1-\rho_N} \frac{1-\rho}{\rho} \frac{\|T^{(k)} - T_0\|}{\|T^{(k)}\|}. \quad (12)$$

In the above simplification, we have assumed that $\frac{\|\phi^{(k+1,n+1)}\|}{\|\phi^{(k+1)}\|} \sim O(1)$ since the perturbation is more sensitive to the neutronics convergence rate ($\rho_N$) and/or the Picard convergence rate ($\rho$). For instance, if the power iteration is used and it is typically very slow, i.e., $\rho_N \sim 1$, then the perturbation can become very large.

In Fourier analysis, we can assume the point-wise uniform convergence such that we drop the norm symbol in Eq. (12):

$$\varepsilon^{(k+1)} \approx \pm\tau \frac{\rho_{PI}}{1-\rho_{PI}} \frac{1-\rho}{\rho} \frac{T^{(k)} - T_0}{T^{(k)}}. \quad (13)$$

Substituting Eq. (13) into (7), we have

$$T^{(k+1)} = T_m + A\Sigma_f(T^{(k)})\phi^{(k+1)}\left[1 \pm \tau C \frac{T^{(k)} - T_0}{T^{(k)}}\right], \quad (14)$$

where the constant $C$ is defined as

$$C = \frac{\rho_N}{1-\rho_N} \frac{1-\rho}{\rho}. \quad (15)$$

## LINEARIZATION

As in [1], we define the following linearized variables:

$$\psi(x,\mu) = \psi_0(x,\mu) + \varepsilon\psi_1(x,\mu), \quad (16a)$$
$$\phi(x) = \phi_0(x) + \varepsilon\phi_1(x), \quad (16b)$$
$$k_{eff} = k_{eff,0}, \quad (16c)$$
$$T(x) = T_0 + \varepsilon T_1(x), \quad (16d)$$
$$\Sigma_i(T) = \Sigma_{i0} + \Sigma_{i1}(T - T_0)$$
$$= \Sigma_{i0} + \varepsilon\Sigma_{i1}T_1(x), \quad i = t,s,f,a \quad (16e)$$

where $\Sigma_{i1}$ is the linear temperature coefficient for each type of cross section.

Substituting Eqs. (16a)–16(c) into Eqs. (6) and (14), after some algebra we obtain the following linearized equations. The subscript "1" in the linearized terms has been dropped.

$$\mu \frac{\partial \psi^{(k+1)}(x,\mu)}{\partial x} + \Sigma_{t0}\psi^{(k+1)}(x,\mu)$$
$$= \frac{1}{2}\Sigma_{t0}\phi^{(k+1)}(x) - \frac{1}{2}\Sigma_{t0}\gamma T^{(k)}(x), \quad (17)$$

$$T^{(k+1)}(x) = A\Sigma_{f0}\phi^{(k+1,n+1)}(x)$$

$$+A\Sigma_{f0}\phi_0\left(\pm\tau\frac{C}{T_0}+\frac{\Sigma_{f1}}{\Sigma_{f0}}\right)T^{(k)}(x), \quad (18)$$

where

$$\gamma = (1-c_0)\left(\frac{\Sigma_{a1}}{\Sigma_{a0}}-\frac{\Sigma_{f1}}{\Sigma_{f0}}\right)\phi_0, \quad (19a)$$

with

$$\Sigma_{a1} = \Sigma_{t1} - \Sigma_{s1}, \quad (19b)$$

$$c_0 = \frac{\Sigma_{s0}}{\Sigma_{t0}}. \quad (19c)$$

**FOURIER ANALYSIS**

As in [1], we introduce the inverse Fourier transforms:

$$\phi^{(k)}(x) = \int_{-\infty}^{+\infty} a^{(k)}(\xi)e^{i\Sigma_{t0}\xi x}d\xi, \quad (20a)$$

$$\psi^{(k)}(x,\mu) = \int_{-\infty}^{+\infty} b^{(k)}(\xi,\mu)e^{i\Sigma_{t0}\xi x}d\xi, \quad (20b)$$

$$T^{(k)}(x) = \int_{-\infty}^{+\infty} c^{(k)}(\xi)e^{i\Sigma_{t0}\xi x}d\xi. \quad (20c)$$

The solutions are required to satisfy the boundary conditions. The discrete Fourier error mode $\xi$ for the reflective boundary conditions are given below.

$$\xi = \frac{\pi}{\Sigma_{t0}L}j, \quad j = \pm 1, \pm 2, \ldots \quad (21)$$

where $L$ is the reactor core height (or the active fuel length), or equivalently the slab thickness in our model problem.

By substituting Eqs. (20a)–(20c) into Eq. (17), after some algebra we obtain

$$a^{(k+1)}(\xi) = -\frac{\gamma\rho_{\text{PI}}(\xi)}{1-\rho_{\text{PI}}(\xi)}c^{(k)}(\xi), \quad (22)$$

where

$$\rho_{\text{PI}}(\xi) = \frac{\tan^{-1}(\xi)}{\xi}. \quad (23)$$

Note that $\rho_{\text{PI}}$ is the spectral radius function for the standard power iteration (PI) method.

Substituting Eqs. (20a) and (20c) into (18), we obtain

$$c^{(k+1)}(\xi) = A\Sigma_{f0}a^{(k+1)}(\xi) + A\Sigma_{f0}\phi_0\left(\pm\tau\frac{C}{T_0}+\frac{\Sigma_{f1}}{\Sigma_{f0}}\right)c^{(k)}(\xi) \quad (24)$$

Substituting Eq. (22) into (24), we have

$$c^{(k+1)}(\xi) = \begin{bmatrix} A\Sigma_{f0}\phi_0\left(\pm\tau\frac{C}{T_0}+\frac{\Sigma_{f1}}{\Sigma_{f0}}\right) \\ -A\Sigma_{f0}\frac{\gamma\rho_{\text{PI}}(\xi)}{1-\rho_{\text{PI}}(\xi)} \end{bmatrix} c^{(k)}(\xi). \quad (25)$$

Then the spectral radius function of the inexact Picard iteration for the coupled N/TH problem is given as

$$\varrho(\xi) = A\Sigma_{f0}\phi_0\left(\pm\tau\frac{C}{T_0}+\frac{\Sigma_{f1}}{\Sigma_{f0}}\right) - A\Sigma_{f0}\frac{\gamma\rho_{\text{PI}}(\xi)}{1-\rho_{\text{PI}}(\xi)}. \quad (26)$$

Substituting Eqs. (3a) and (19a) into (26), we obtain

$$\varrho(\xi) = \pi r_{fo}^2\kappa R_t\Sigma_{f0}\phi_0\left[\begin{array}{c}\pm\tau\frac{C}{T_0}+\frac{\Sigma_{f1}}{\Sigma_{f0}}\\-\left(\frac{\Sigma_{a1}}{\Sigma_{a0}}-\frac{\Sigma_{f1}}{\Sigma_{f0}}\right)\frac{1-c_0}{1-\rho_{\text{PI}}(\xi)}\rho_{\text{PI}}(\xi)\end{array}\right]. \quad (27)$$

Noting $q_0' = \pi r_{fo}^2\kappa\Sigma_{f0}\phi_0$, Eq. (27) can be written as

$$\varrho(\xi) = q_0'R_t\left[\pm\tau\frac{C}{T_0}+\frac{\Sigma_{f1}}{\Sigma_{f0}}-\left(\frac{\Sigma_{a1}}{\Sigma_{a0}}-\frac{\Sigma_{f1}}{\Sigma_{f0}}\right)\frac{1-c_0}{1-\rho_{\text{PI}}(\xi)}\rho_{\text{PI}}(\xi)\right]. \quad (28)$$

Given that $q_0'R_t = T_0 - T_m$, Eq. (28) can be rewritten as

$$\varrho(\xi) = (T_0-T_m)\left[\pm\tau\frac{C}{T_0}+\frac{\Sigma_{f1}}{\Sigma_{f0}}-\left(\frac{\Sigma_{a1}}{\Sigma_{a0}}-\frac{\Sigma_{f1}}{\Sigma_{f0}}\right)\frac{1-c_0}{1-\rho_{\text{PI}}(\xi)}\rho_{\text{PI}}(\xi)\right], \quad (29)$$

where the constant $C$ is defined by Eq. (15).

Finally, the spectral radius of Picard iteration for the coupled N/TH nonlinear system is given as

$$\rho = \max_\xi|\varrho(\xi)|. \quad (30)$$

As compared to the standard Picard iteration that fully converges the neutronics solution during iteration [1], the inexact Picard iteration has an extra perturbation term, $\frac{C}{T_0}$, where $C = \frac{\rho_N}{1-\rho_N}\frac{1-\rho}{\rho}$. If the perturbation is too large (e.g., $\rho_N\sim 1$), then a sufficiently small forcing parameter $\tau$ must be used to ensure the stability of the inexact Picard iteration, that is, the neutronics equation should be "over solved" with sufficient accuracy during iteration.

**RESULTS**

To verify the theoretical analysis, we compare the FA predictions with numerical results based on a 1-D model problem, which is the homogeneous slab with the reflective boundary on both sides. The Gauss-Legendre S12 quadrature set is used for angular discretization and the Diamond Difference (DD) method is employed for spatial discretization. Note that the angular quadrature and the mesh size used are sufficiently fine to minimize the numerical errors. A simple heat balance model is used to calculate the fuel and coolant temperatures at each axial cell. For this case, the core height $L = 150$ cm. The cross section data for a typical PWR fresh fuel pin is summarized in Table I.

TABLE I. Typical PWR Data

| $\Sigma_{t0}$ (cm$^{-1}$) | $\nu\Sigma_{fo}$ (cm$^{-1}$) | $c_o$ | $\Sigma_{f1}/\Sigma_{f0}$ (K$^{-1}$) | $\Sigma_{a1}/\Sigma_{a0}$ (K$^{-1}$) |
|---|---|---|---|---|
| 0.534 | 0.0255 | 0.96 | $-2.59\times 10^{-5}$ | $9.63\times 10^{-6}$ |

The temperature coefficient of the absorption cross section is positive, and it is negative for the fission cross section. For such problems, the convergence of the inexact Picard iteration is determined by the smallest error mode $\xi = \pi/(\Sigma_{t0}L)$, and the spectral radius is given as

$$\rho = (T_0-T_m)\left[\tau\frac{C}{T_0}-\frac{\Sigma_{f1}}{\Sigma_{f0}}+\left(\frac{\Sigma_{a1}}{\Sigma_{a0}}-\frac{\Sigma_{f1}}{\Sigma_{f0}}\right)\frac{1-c_0}{1-\rho_{\text{PI}}(\xi)}\rho_{\text{PI}}(\xi)\right]. \quad (31)$$

When the power iteration is used to obtain the neutronics solution, the FA predicts the spectral radius of the neutronics iteration, $\rho_N = \rho_{\text{PI}} = \frac{\tan^{-1}(\pi/(\Sigma_{t0}L))}{\pi/(\Sigma_{t0}L)} = 0.99949$, which is almost the same as the numerical result. It should be noted that if a different iterative method (other than the power iteration) is used for the neutronics solution, then $\rho_N$ is generally different with $\rho_{PI}$. The FA calculated spectral radius of Picard iteration, $\rho = 0.745$, while it is 0.737 from the numerical calculation. The FA predicts that Picard iteration becomes unstable when $\tau > 0.0012$ and the numerical result is about 0.0015.

The perturbation $\frac{C}{T_0}$, where $C = \frac{\rho_N}{1-\rho_N}\frac{1-\rho}{\rho}$, is sensitive to the convergence rate of the neutronics iteration. For this model problem, $C = \frac{0.99949}{1-0.99949}\frac{1-0.745}{0.745} = 678.8$, and then $\frac{C}{T_0} = \frac{678.8}{850} = 0.79$, while the sum of other terms in Eq. (31) is only 0.0027. Note that $T_0 - T_m = 275$ K. Therefore, it requires that $\tau$ be sufficiently small to minimize the impact of the perturbation on the stability.

We next use a nonlinear diffusion acceleration scheme, lpCMFD [7], to speed up the power iteration. The calculated spectral radius of the accelerated power iteration, $\rho_N = 0.63$. We recalculate the perturbation: $C = \frac{0.63}{1-0.63}\frac{1-0.745}{0.745} = 0.58$, and then $\frac{C}{T_0} = \frac{0.58}{850} = 6.82 \times 10^{-4}$, which is now smaller than 0.0027. Therefore, we can use a relatively large value for the forcing parameter to reduce over-solving.

Finally, we demonstrate that the adaptive termination criterion can effectively reduce the computational cost when the neutronics iteration is accelerated by lpCMFD. Fig. 1 shows the number of neutronics iterations during each Picard iteration for different values of $\tau$. For this case, $\varepsilon = 10^{-7}$ is used for the fuel temperature tolerance, and it is $10^{-8}$ for the final neutron flux tolerance. The refence case is depicted by the blue curve, in which a fixed flux stopping criterion ($10^{-8}$) is used.

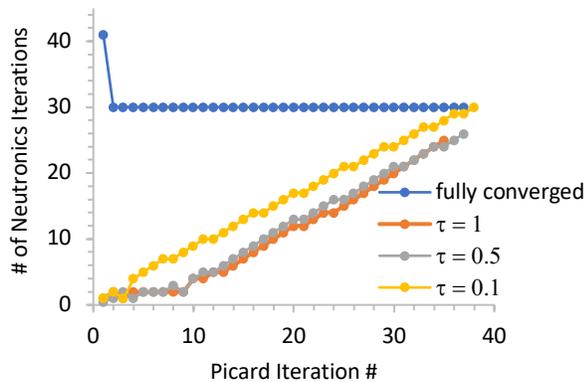

Fig. 1. Adaptive inexact Picard iteration.

As discussed above, $\tau$ can be relatively large due to small perturbation error, and thus it is very effective in mitigating over-solving, especially for the earlier Picard iterations. Overall, the adaptive termination can save almost 50% of neutronics iterations for this case. For the fully converged case, more neutronics iterations are taken at the first Picard iteration because the diffusion acceleration is less effective for the flat flux solution (given the constant initial cross sections and the reflective BCs).

## CONCLUSIONS

We have presented a formal Fourier analysis to characterize the convergence properties of the inexact Picard fixed-point iteration with adaptive termination tolerance for the coupled N/TH calculation. The inexact Picard iteration has been modeled as a perturbed fixed-point iteration to determine the forcing parameter. It is shown that this parameter is sensitive to the convergence rate of the inner neutronics solution. If the neutronics solution converges very slowly (i.e., $\rho_N \sim 1$), then a sufficiently small number must be used to ensure the stability. It means that a certain amount of over-solving cannot be avoided. However, when the neutronics iteration is accelerated with an acceleration scheme (e.g., lpCMFD), the restriction on the forcing parameter can be greatly relaxed and a moderate value (e.g., $\tau \sim 0.5$) can be very effective in improving the computational cost.

If the TH solution is obtained iteratively, the Fourier analysis presented can be extended to determine the proper forcing parameter in the adaptive termination criterion for the TH calculation.

## ACKNOWLEDGEMENT

The second author was supported by the Exascale Computing Project (17-SC-20-SC), a collaborative effort of the U.S. Department of Energy Office of Science and the National Nuclear Security Administration.

## REFERENCES

1. D. WANG, "Stability Analysis of Picard Iteration for Coupled Neutronics/Thermal-Hydraulics Simulations," *Trans. Am. Nucl. Soc.*, **128**, 246–249 (2023).
2. R. S. DEMBO, et al., "Inexact Newton Methods," *SIAM J. Numer. Anal.*, **19**(2), 400–408 (1982).
3. S. C. EISENSTAT and H. F. WALKER, "Choosing the Forcing Terms in an Inexact Newton Method," *SIAM J. Sci. Comput.*, **17**(1), 16–32 (1996).
4. P. BIRKEN, "Termination Criteria for Inexact Fixed-Point Schemes," *Numer. Linear Algebra Appl.*, **22**, 702–716 (2015).
5. J. P. SENECAL, W. JI, "Development of an Efficient Tightly Coupled Method for Multiphysics Reactor Transient Analysis," *Prog. Nucl. Energy*, **103**, 33–44 (2018).
6. J. M. ORTEGA and W. C. RHEINBOLDT, *Iterative Solution of Nonlinear Equations in Several Variables*, SIAM, Philadelphia, PA, USA (2000).
7. D. WANG and S. XIAO, "A Linear Prolongation Approach to Stabilizing CMFD," *Nucl. Sci. Eng.*, **190**, 45–55 (2018).